
\magnification=\magstep 1
\baselineskip=15pt     
\parskip=3pt plus1pt minus.5pt
\font\hd=cmbx10 scaled\magstep1
\overfullrule=0pt
\centerline{\hd A Finiteness Theorem for Elliptic Calabi-Yau Threefolds}
\input cyracc.def
\newfam\cyrfam
\font\tencyr=wncyr10
\def\cyr{\fam\cyrfam\tencyr\cyracc}
\def\TS{\hbox{\cyr Sh}}
\def\num{\global\advance\count10 by 1 \eqno(\the\count10)}

\def\Ptwo{{\bf P^2}}
\def\P{{\bf P}}
\def\C{{\cal C}}

\def\SS{{\cal S}}
\def\X{{\cal X}}
\def\O{{\cal O}}
\def\L{{\cal L}}
\def\T{{\cal T}}

\def\G{{\cal G}}

\def\F{{\cal F}}

\def\E{{\cal E}}

\def\ZZ{{\cal Z}}

\def\Pic{\mathop{\rm Pic}}
\def\Spec{\mathop{\rm Spec}}

\def\boldz{{\bf Z}}

\def\Pone{{\bf P}^1}

\def\Gm{{\bf G}_m}

\def\qz{{\bf Q}/{\bf Z}}

\def\exact#1#2#3{0\rightarrow#1\rightarrow#2\rightarrow#3\rightarrow0}

\def\mapleft#1{\smash{
  \mathop{\longleftarrow}\limits^{#1}}}
\def\mapright#1{\smash{
  \mathop{\longrightarrow}\limits^{#1}}}
\def\mapsw#1{\smash{
  \mathop{\swarrow}\limits^{#1}}}
\def\mapse#1{\smash{
  \mathop{\searrow}\limits^{#1}}}
\def\mapdown#1{\Big\downarrow
   \rlap{$\vcenter{\hbox{$\scriptstyle#1$}}$}}
\medskip
\centerline{\it Mark Gross\footnote{*}{
Research at MSRI supported in part
by NSF grant \#DMS 9022140.}}
\medskip
\centerline{May, 1993}
\medskip
\centerline{Mathematical Sciences Research Institute}
\centerline{1000 Centennial Drive}
\centerline{Berkeley, CA 94720}
\centerline{mgross@msri.org}
\bigskip
\bigskip
{\hd \S 0. Introduction}

For the purposes of this paper, we define a
{\it Calabi-Yau threefold} to be an algebraic threefold $X$ over
the field of complex numbers which is birationally
equivalent to a threefold $Y$ with ${\bf Q}$-factorial
terminal singularities, $K_Y=0$, and $\chi(\O_Y)=h^1(Y,\O_Y)=h^2(Y,\O_Y)=0$.
We say that
$Y$ is a minimal Calabi-Yau threefold.

The analagous class for surfaces are the K3 surfaces.
All K3 surfaces are homeomorphic: there is one underlying topological
type. On the other hand, there are a large number of topological types
of  minimal Calabi-Yau threefolds, but it is an open question of
whether there are a finite
number of such types.
A stronger question would be to ask whether there are a finite number of
families of algebraic minimal Calabi-Yau threefolds. This is definitely
not true for K3 surfaces: there are a countably infinite number of algebraic
families. Up to birational equivalence, we answer this question for those
Calabi-Yaus which possess an elliptic fibration.

Our main theorem is:

\proclaim Theorem 0.1. There exists a finite number of triples
$(\X_i,\SS_i,\T_i)$ of quasi-projective varieties with maps
$$\matrix{\X_i&&\cr
\mapdown{f_i}&\mapse{\pi_i}&\cr
\SS_i&\mapright{g_i}&\T_i\cr}$$
where $\pi_i$ is smooth and proper
with each fibre a Calabi-Yau threefold,
$f_i$ proper with generic fibre an elliptic curve, and $g_i$
smooth and proper
with each fibre a rational surface,
such that for any elliptic
fibration $X\rightarrow S$ with $X$ Calabi-Yau and $S$ rational
there exists a $t\in\T_i$ for some $i$ such
that there are birational maps $X\cdot\cdot\rightarrow (\X_i)_t$,
$S\cdot\cdot\rightarrow (\SS_i)_t$ with the following diagram commutative:
$$\matrix{X&\cdot\cdot\rightarrow&(\X_i)_t\cr
\mapdown{}&&\mapdown{}\cr
S&\cdot\cdot\rightarrow&(\SS_i)_t\cr}$$

Let us make several remarks before indicating the idea of the proof.
First, in this theorem, we only consider elliptic fibrations with rational
bases. If $X$ is Calabi-Yau and $f:X\rightarrow S$ an elliptic fibration,
then $S$ is either rational or birational to an Enriques surface. In the
latter case, $X$ has a particularly simple structure.
In particular, it is the quotient of a three-fold $Y$ with
$\kappa(Y)=0$ and $h^1(\O_Y)>0$. (See Proposition 2.10).

Secondly, note that by [12], Theorem 1.2,
we can thus obtain a finite number of
families of minimal Calabi-Yau threefolds $\X'_i\rightarrow\T_i$. However,
because of the non-uniqueness of minimal models for threefolds,
not every minimal elliptic Calabi-Yau threefold will be isomorphic
to a $(\X'_i)_t$ for some $t\in\T_i$, but merely birational to such.
It is not known if a Calabi-Yau threefold always has a finite number
of minimal models, even up to automorphism.
Thus this result does not imply that there are only a finite number of
topological types of even non-singular minimal elliptic Calabi-Yau threefolds.
Nevertheless, any such threefold will be related by a series of flops to
a finite number of possible topological types.

Thirdly, this result is a much stronger one than proved in [10], i.e.
that there exists a finite number of possible Euler characteristics of certain
types of elliptic Calabi-Yau threefolds. Our theorem says that there exists
a complete, finite classification. One hopes that any
minimal Calabi-Yau threefold with sufficiently large Picard
number is elliptic. With our results, proving such a conjecture
would then show that there are only a finite number of types
of minimal Calabi-Yau threefolds with Picard number greater than some fixed
number. This hopefully gives some suggestion that there are only a finite
number of types of algebraic Calabi-Yau threefolds.

If one were to ask how many families of elliptic Calabi-Yaus exist, the
proof of this theorem gives no hint, other than to suggest that it
is a very large number. Making an uninformed guess, I would conservatively
expect thousands of families.

The proof of this theorem relies heavily on results about elliptic threefolds.
Results of [3] and [7] tell us how to go about classifying elliptic
three-folds; we simply have to apply them carefully. The proof also uses
minimal model theory for threefolds, especially as developed for
elliptic threefolds in [4] and [6]. The main application of these
techniques is

\proclaim Theorem 0.2. Let $f:X\rightarrow S$ be an elliptic fibration with $X$
and $S$ projective and
$X$ having a minimal model with trivial canonical class. Then
there exists a
birationally equivalent equidimensional
fibration $\bar f:\bar X\rightarrow \bar S$
where $\bar X$ is
a minimal model of $X$ and $\bar S$ is a projective surface with only DuVal
singularities.

Grassi in [4] has proven
this theorem in the case that $f$ has no multiple fibres. In fact, in
that case $\bar S$ can be found to be non-singular.

The second key ingredient is the calculation of the Tate-Shafarevich group
for elliptic fibrations in [3] and [7]. The most important fact is that
the Tate-Shafarevich group is finite for Calabi-Yau elliptic fibrations,
something which fails for K3 surfaces. After these
two observations, most of the proof is simply careful bookkeeping to make
sure that
we have not missed any threefolds.

I would like to thank M. Reid for suggesting this problem to me, and I.
Dolgachev and
A. Grassi for many useful discussions.

{\hd \S 1. Generalities on Elliptic Threefolds.}

We recall some basic definitions and state some general results about
elliptic threefolds which will be necessary.
We note that all varieties in this paper
are algebraic varieties of finite type over the complex
numbers.

\proclaim Definition 1.1. A projective morphism $f:X\rightarrow S$ is
called an {\it elliptic fibration} if its generic fibre $E$ is
a regular curve of genus one and all fibres are geometrically connected.
$f:X\rightarrow S$ is called a model for $E$.
The closed subset
$$\Sigma_{red}=\Sigma_{red}(f)=\{s\in S|\hbox{$X_s$ is not regular}\}$$
is called the {\it (reduced) discriminant locus}. If $t\in S$ is a codimension
one point of $S$, then
the {\it fibre type of $t$} is the Kodaira fibre type ($_mI_n,I_n^*,II,II^*,
III,III^*,IV,IV^*$) of the central fibre of a relatively minimal model (or
N\'eron model) of $X(\bar t)=X\times_S \Spec \O_{S,\bar t}$. ($\O_{S,\bar t}$
is the strict henselization of the local ring $\O_{S,t}$.)
A {\it collision} is a singular point of $\Sigma$. The closed subset
$$\Sigma^m=\{s\in S| \hbox{$f$ is not smooth at any $x\in f^{-1}(s)$}\}$$
is called the {\it multiple locus} of $f$. A fibre over a point
$s\in\Sigma^m$ is called {\it multiple}. A fibre is called an
{\it isolated multiple fibre} if it is over a zero-dimensional
component of $\Sigma^m$.
A {\it section} (resp. a {\it rational
section}) of $f$ is a closed subscheme $Y$ of $X$ for which the
restriction of $f$ to $Y$ is an isomorphism (resp. a birational morphism).

Let $f:X\rightarrow S$ be an elliptic fibration between two non-singular
complex varieties, which is smooth off of a simple normal crossings
divisor on $S$. Then by [11], Theorem 20, $f_*\omega_{X/S}$ is an
invertible sheaf. This defines a divisor on $S$, which we denote by
$\Delta_{X/S}$ (or $\Delta$ if no confusion will arise). If $f:X\rightarrow
S$ is then just an elliptic fibration with $X$ and $S$ having singularities
in codimension 3 and 2 respectively, then there exists an open subset
$S_0\subseteq S$ such that $S-S_0$ is codimension 2 in $S$ and if $f_0$ is the
induced morphism $X_0=
X\times_S S_0\rightarrow S_0$, then $X_0$ and $S_0$ are non-singular
and $\Sigma_{red}(f_0)$ has simple normal crossings. Then $\Delta_{X_0/S_0}$
is a Cartier divisor on $S_0$, which extends to a Weil divisor on $S$, which we
denote by $\Delta_{X/S}$ again.

Let $\{M_i\}$ be the codimension one components of $\Sigma^m$, with
fibre type ${}_{m_i}I_a$. We set
$$\Lambda_{X/S}=\Delta_{X/S}+\sum {m_i-1\over m_i} M_i.$$
We then have

\proclaim Lemma 1.2. Let $f:X\rightarrow S$ be an elliptic fibration
with $\dim S=2$.
Then there exists a blow-up $S'\rightarrow S$ and a birationally equivalent
fibration $f':X'\rightarrow S'$
such that
\item{1)} $f'$ is flat and relatively minimal,
$X'$ has only ${\bf Q}$-factorial terminal singularities,
$S'$ is regular, and the
reduced discriminant locus $\Sigma_{red}$ of $f'$ has simple
normal crossings.
\item{2)}
The modular function $J$
gives
a morphism $J:S'\rightarrow \Pone$, and $\Delta=\Delta_{X'/S'}$ is
a divisor on $S'$
such that
$$12\Delta\sim J_{\infty}+\sum 12a_iD_i$$
where $J_{\infty}$ is the fibre of $J$ at $\infty\in\Pone$, and $D_i$
are the components of $\Sigma_{red}$
with $12a_i=0,2,3,4,6,8,9$ or $10$ depending
on whether $D_i$ is of fibre-type ${}_mI_a,II,III,IV,I_a^*,IV^*,III^*$ or
$II^*$
respectively.
\item{3)} $$K_{X'}=f'^*(K_{S'}+\Lambda_{X'/S'}).$$
Since $f'^*M_i$ is a divisor with multiplicity
$m_i$, this gives a well-defined ${\bf Q}$-cartier Weil divisor on $X'$.

Proof: By [20], Theorem A.1, there exists a blowing up $S'\rightarrow S$
and a model $X'\rightarrow S'$ satisfying 1). (We can always assume that
$\Sigma_{red}$ has simple normal crossings by blowing up $S$ first until this
is achieved.) 2) then follows from [11], Theorem 20.
3) follows from [19], Theorem 0.1. $\bullet$

Next, let $\phi:S\rightarrow S'$ be a birational morphism from a smooth
surface $S$ to a possibly singular surface $S'$, with normal
singularities, with exceptional curves $E_1,\ldots,E_n\subseteq S$. Recall
that if $D$ is a ${\bf Q}$-Weil divisor on $S'$, then $\phi^*D$ is the
unique ${\bf Q}$-divisor on $S$ of the form $\phi^{-1}(D)+\sum a_iE_i$
for some $a_i$ such that $\phi^*D.E_i=0$
for all $n$.  This allows us to define an intersection product on $S'$ by
$C.D=\phi^*C.\phi^*D$.

The following theorem is an application of Mori's minimal model algorithm,
and will be crucial for the classification of singularities of the base.

\proclaim Lemma 1.3. Let $f:X\rightarrow S$ be a relatively minimal elliptic
fibration with $K_X=f^*(K_S+\Lambda_{X/S})$, $S$ a non-singular surface.
Let $\phi:S\rightarrow S'$ be a birational morphism
with $S'$ normal, exceptional curves
$E_1,\ldots,E_n$, and $K_S+\Lambda_{X/S}=\phi^*(K_{S'}
+\Lambda_{X/S'})+\sum_{i=1}^n a_iE_i$ with $a_i>0$ for all $i$. Then if
$f':X'\rightarrow S'$ is a relatively minimal model of $\phi\circ f:
X\rightarrow S'$, then $K_{X'}=f'^*(K_{S'}+\Lambda_{X'/S'})$. (Note that
$\Lambda_{X'/S'}=\Lambda_{X/S'}$.) Furthermore,
if $f$ is equidimensional, then $f'$ is equidimensional.

Proof: [6], Theorem 2.5. $\bullet$

\proclaim Definition 1.4. Let $f:X\rightarrow S$ be an elliptic fibration.
An elliptic fibration $j:J\rightarrow S$ is called the jacobian of $f$ if
the generic fibre of $j$ is the jacobian of the generic fibre
of $f$. The jacobian of $f$ is
thus defined up to birational equivalence.

Any jacobian fibration over a non-singular variety $S$
is birationally equivalent to a Weierstrass model over $S$ ([3], Prop. 2.4).
We review the definition of a Weierstrass model. See
[1,18,19] for details.
Let $\L$ a line bundle on a scheme $S$, $a\in H^0(S,\L^{\otimes 4})$,
and $b\in H^0(S,\L^{\otimes 6})$ such that $4a^3+27b^2$ is a non-zero
section of $\L^{\otimes 12}$. Let $\P=\P(\O_S\oplus \L^{\otimes-2}
\oplus \L^{\otimes-3})$, $\pi:\P\rightarrow S$ be the natural projection,
and $\O_{\P}(1)$  the tautological line bundle on $\P$.
We define the scheme $W(\L,a,b)$ as a closed subscheme of $\P$
given by the equation $Y^2Z=X^3+aXZ^2+bZ^3$
where  $X$, $Y$ and $Z$
are given by the sections of
$\O_{\bf P}(1)\otimes \L^{\otimes 2}$,
$\O_{\bf P}(1)\otimes \L^{\otimes 3}$, and
$\O_{\bf P}(1)$ which correspond to the natural injections
of $\L^{\otimes-2}$, $\L^{\otimes-3}$ and $\O_S$ into
$\pi_*\O_{\P}(1)=
\O_S\oplus\L^{\otimes -2}\oplus\L^{\otimes -3}$, respectively.

The structure morphism $f:W(\L,a,b)\rightarrow S$ is a flat elliptic fibration,
called a {\it Weierstrass fibration}. It has a section
$\sigma:S\rightarrow W(\L,a,b)$ defined by the $S$-point
$(X,Y,Z)=(0,1,0)$. We also see that $\sigma(S)$ lies in the smooth
locus of $W(\L,a,b)$ if $S$ is regular.
We will call this section the {\it section at infinity}. It is easy to
see that $K_{W(\L,a,b)}=f^*(K_S\otimes\L)$.

A Weierstrass fibration $W(\L,a,b)\rightarrow S$ is called minimal if there
is no effective divisor $D$ on $S$ such that $div(a)\ge 4D$, $div(b)\ge 6D$.
Any Weierstrass model is birationally equivalent to a minimal Weierstrass
model.

The reduced discriminant locus $\Sigma_{red}$
of $W(\L,a,b)\rightarrow S$ is equal to the
support of the Cartier divisor defined by the section $\Sigma$
of $\L^{\otimes 12}$ given by $4a^3+27b^2$.
This gives the discriminant locus a scheme structure.
If $j:J\rightarrow S$ is a relatively minimal jacobian fibration and
$w:W(\L,a,b)\rightarrow S$ a birationally equivalent minimal
Weierstrass model, then $\Sigma_{red}(j)=supp(\Sigma(w))$. Thus, when
dealing with a relatively minimal jacobian fibration, we can set
$\Sigma(j):=\Sigma(w)$, and $\Sigma(j)\sim 12\Delta_{J/S}$.

Let $W(\L,a,b)\rightarrow S$ be a Weierstrass
model with $\dim S=2$.
Then there exists a blowing-up $S'\rightarrow S$ with $S'$ regular, a
Weierstrass model $W(\L',a',b')\rightarrow S'$ birational
to $W(\L,a,b)\rightarrow S$, and a resolution of singularities $X'\rightarrow
W(\L',a',b')$ with
$X'$ flat over $S'$.
Indeed, Miranda [15] has given an explicit algorithm for finding
such a resolution,
first by blowing up the base surface $S$ until the reduced
discriminant
locus $\Sigma_{red}$ has simple normal crossings, and
continuing further so that only one of a small list of possible
collisions between components of $\Sigma$ can occur, namely
the following possibilities:
$I_{M_1}+I_{M_2},
I_{M_1}+I_{M_2}^*,
II+IV,
II+I_0^*,
II+IV^*,
IV+I_0^*,
III+I_0^*$. Thus we have

\proclaim Definition 1.5. A {\it Miranda elliptic fibration} is
an elliptic fibration $f:X\rightarrow S$ such that
\item{a)} $X$ and $S$ are regular and $f$ is flat and has a section;
\item{b)} the reduced discriminant locus $\Sigma_{red}$
has simple normal crossings;
\item{c)} All collisions are of type $I_{M_1}+I_{M_2}$, $I_{M_1}+I_{M_2}^*$,
$II+IV$, $II+I_0^*$, $II+IV^*$, $IV+I_0^*$ or $III+I_0^*$.

A few comments about passing from an elliptic three-fold to its jacobian.
If $f:X\rightarrow S$ and $j:J\rightarrow S$ are relatively minimal, then
outside a codimension 2 subset and $\Sigma^m$, the two fibrations are
locally (i.e. in the complex or \'etale topologies on $S$) isomorphic. This
follows from the proof of [3], Prop. 2.17. Fibres over collision
points may change; in particular, there may be flat relatively minimal
models $f:X\rightarrow S$ such that there is no flat relatively minimal
model for its jacobian over $S$. Fortunately, the canonical bundle formula
remains valid for a relatively minimal model of the
jacobian even though it may not be flat, assuming the discriminant locus
has simple normal crossings.

\proclaim Lemma 1.6. Let $f:X\rightarrow S$ be a relatively minimal
elliptic fibration  with $\Sigma_{red}(f)$ simple normal crossings. Then
a relatively minimal model $j:J\rightarrow S$ of the jacobian of $f$
obeys
$$K_{J}=j^*(K_{S}+\Delta_{J/S}).$$
Furthermore, $12\Delta_{J/S}=12\Delta_{X/S}$.

Proof. Since $\Sigma_{red}(j)\subseteq \Sigma_{red}(f)$, we see that
$\Sigma_{red}(j)$ must have simple normal crossings. Thus by [6],
Theorem 1.14, the formula for the canonical class of $J$ holds.

For the second statement, it is clear that the $J$-morphism for $f$ and $j$
coincide, since $f$ and $j$ are locally isomorphic away from the multiple
fibre locus. Furthermore, in codimension one the only multiple
fibres possible are of type ${}_mI_a$. By [11], Theorem 20,
we have
$$12\Delta_{X/S}=J_{\infty}+\sum a_iD_i,$$
where the $a_i$ and $D_i$ are as in item 2) of Lemma 1.2. The same
formula holds for $12\Delta_{J/S}$, and by the above observations,
$J_{\infty}$ and the $a_i$ and $D_i$ coincide for both fibrations. $\bullet$

{\hd \S 2. Proof of Theorem 0.2 and Related Results.}

{\it Proof of Theorem 0.2}:
Let $f:X\rightarrow S$ be an elliptic fibration with $X$ birational
to a threefold with trivial canonical class. Replace $f:X\rightarrow S$
with the birationally equivalent fibration given in Lemma 1.2.

Following [4], in the proof of Lemma 1.4,
one can then find a unique contraction map
$\phi:S\rightarrow \bar S$ such that the Zariski decomposition of
the ${\bf Q}$-divisor
$K_{S}+\Lambda_{X/S}$ is $\phi^*(K_{\bar S}+\Lambda_{X/\bar S})+\sum
c_iE_i$, the $E_i$ the exceptional
curves of $\phi$, with $K_{\bar S}+\Lambda_{X/\bar S}$
nef, the $c_i$ positive rational numbers, and
$\bar S$ has log-terminal singularities.
Then Lemma 1.3  shows that if one takes a relatively minimal model
$\bar f:\bar X\rightarrow \bar S$ of $X\rightarrow \bar S$,
then $\bar X$ is in fact minimal, with $K_{\bar X}=\bar f^*(K_{\bar
S}+\Lambda_{
\bar X/\bar S})$,
and thus in our case $K_{\bar S}+\Lambda_{\bar X/\bar S}=0$.
Furthermore $\bar f$ is equidimensional.
In the proof of [21] Theorem 3.1, Oguiso
shows that  if $K_{\bar X}=0$, then $\bar f$
cannot have multiple fibres in codimension one; hence it has only isolated
multiple fibres. Thus the one-dimensional components $M_i$ of $\Sigma^m(f)$
are all exceptional
for $\phi$, and we can write $0=K_{\bar X}=
\bar f^*(K_{\bar S}+\bar\Delta)$,
where $\bar\Delta=\Delta_{X/\bar S}=\Delta_{\bar X/\bar S}$.

We need to show that $\bar S$ has only canonical
(i.e. DuVal) singularities. Let $P\in \bar S$ be a singular point, and
replace $S$ and $\bar S$ with analytic germs around $\phi^{-1}(P)$ and $P$.
We have $\phi:S\rightarrow \bar S$, with
$\phi^{-1}(P)=\bigcup_{i=1}^r E_i$, each $E_i$ a $\Pone$, as $P$ is necessarily
a log-terminal, hence rational, singularity. As $S$ and $\bar S$ are germs,
$$\Pic S=H^1(S,\O_S^*)=H^2(S,\boldz)=\boldz^r,$$
generated by divisors $D_i$, $i=1,\ldots,r$, with $D_i.E_j=\delta_{ij}$.
We have
$$\eqalign{Cl(\bar S)&=\Pic (\bar S-P)\cr
&=\Pic (S-\phi^{-1}(P))\cr
&={{\boldz^r}\over (E_1,\ldots,E_r)},\cr}$$
where $(E_1,\ldots,E_r)$ is the subgroup of $\Pic S$ generated by the
$E_i$.
The intersection matrix of the $E_i$'s is negative definite, so we can
write any element $D\in \Pic(S)$ as a linear combination of the $E_i$ with
rational coefficients. $\phi_*:\Pic(S)\rightarrow Cl(\bar S)$
is defined to be
the projection $\boldz^r\rightarrow \boldz^r/(E_1,\ldots,E_r)$,
so $\phi_*D=0$ in $Cl(\bar S)$ if the coefficients of the $E_i$'s
in $D$ are integral. It is clear that $\bar \Delta=\phi_*\Delta_{X/S}$.
Thus we can write
$$K_{S}+\Delta_{X/S}+\sum{m_i-1\over m_i}E_i
=\sum_{i=1}^r\left(d_i+a_i+{m_i-1\over m_i}\right) E_i+b_i D_i.$$
Here, $d_i$ is the discrepancy of $E_i$ (i.e. $K_{S}=\sum d_iE_i$), $12a_i$ is
the appropriate coefficient  of $E_i$ in $12\Delta_{X/S}$ in 2) of Lemma 1.2,
$a_i$ depending on the
fibre type of $E_i$, $0\le a_i <1$;
$m_i$ is the multiplicity of fibres over $E_i$.
Notice
that if $a_i>0$, then $m_i=1$, since in codimension one
there are only multiple fibres of type
${}_mI_a$. Thus
$$0\le a_i+{m_i-1\over m_i}<1. \leqno{(2.1)}$$
Finally $b_i$ is the suitable coefficient for those components of
$\Sigma_{red}(f)$ appearing in $\Delta_{X/S}$
which are not one of the $E_i$'s but intersect them (this includes
$J_{\infty}$,
which can be replaced by another, linearly equivalent,
fibre of $J$ which does not contain any
$E_i$).
We have $0\le b_i$, and
we can write $\sum b_i D_i=\sum b_i'E_i$, with $b_i'\le 0$, by negative
definiteness of the intersection matrix of the $E_i$'s. In fact, since $\sum
b_iD_i$ is nef, $b_i'<0$ for all $i$ unless $\sum b_iD_i=0$, i.e. $b_i=0$
for all $i$. (See, for example, [13], 2.19.3).

By uniqueness of the Zariski decomposition, we must have
$$c_i=d_i+a_i+{m_i-1\over m_i}+b_i'.$$

Furthermore, we have
$$0=K_{\bar S}+\bar\Delta=\phi_*(K_S+\Delta_{X/S}),$$
which means that
$$d_i+a_i+b_i'\in\boldz\quad\forall i\leqno{(2.2)}$$
Furthermore,
$$d_i+a_i+{m_i-1\over m_i}+b_i'>0 \quad \forall i\leqno{(2.3)}$$
since $c_i>0$ for all $i$.

Now suppose that $d_i<0$ for some $i$. Then $d_i+b_i'<0$,
so (2.2) and $a_i<1$
imply that $d_i+b_i'+a_i\le 0,$ with equality only possible
if $a_i>0$, and hence $m_i=1$, contradicting (2.3). If
$a_i=0$, then $d_i+b_i'\le -1$, again contradicting (2.3).
Thus $d_i\ge 0$ for all $i$,
and $\bar S$ has a canonical, hence DuVal, singularity. This
proves Theorem 0.2.

Note that if $m_i=1$ for all $i$, then if $d_i=0$, $b_i'+a_i$ cannot be
greater than zero, so $d_i>0$ for all $i$, and hence $\bar S$ has
terminal singularities, i.e., $\bar S$ is smooth, as was shown in
[4], Cor. 3.3, using results in [5]. $\bullet$

\proclaim Corollary 2.1. In the notation of the proof of 0.2,
let $E_1,\ldots,E_r$ be the set of curves in $S$ mapping to a singular point
$P$
of $\bar S$. Then $\bigcup E_i$ is contained in a fibre of the $J$-morphism.
Furthermore, if a relatively minimal model for the jacobian $j:J\rightarrow
S$
has a singular fibre over any point
in $\bigcup E_i$, then $\bigcup E_i$ is contained in $J_{\infty}$.

Proof: Assume that $\bigcup E_i$ is not contained in a fibre of the
$J$-morphism. Then one of the $b_i$ would have to be non-zero, as at least
one of the $E_i$ intersects $J_{\infty}$. Let $k$ be chosen so that $d_k=0$.
If $b_k'<0$, then $d_k+b_k'<0$ and (2.2) then
implies that $a_k>0$ so that $m_k=1$, contradicting (2.2) and
(2.3). Thus $b_k'=0$.
This then
implies ([13, 2.19.3]) that $b_i=0$ for all $i$.

For the second statement, if $j:J\rightarrow S$ had a
singular fibre over some point
of $\bigcup E_i$ and $\bigcup E_i$ were not contained in $J_{\infty}$,
then as $12\Delta_{X/S}=12\Delta_{J/S}$ (Lemma 1.6)
either $b_i\not=0$ for some $i$, which we have already
seen is not possible, or else $a_i>0$ for some $i$. But all the $d_i$ are
integers, and $a_i$ is only an integer if $a_i=0$, so we obtain a contradiction
by (2.2).
$\bullet$

\proclaim Proposition 2.2. Let $f:X\rightarrow S$ be an elliptic fibration
with $X$ a Calabi-Yau threefold and $S$ rational, and let $j:J\rightarrow
S$ be its jacobian fibration. Then $J$ is a Calabi-Yau threefold,
and there is a minimal model $J'$ of $J$ with a fibration $j':J'\rightarrow
S'$ along with a morphism $S'\rightarrow \bar S$, where $\bar f:
\bar X\rightarrow \bar S$ is the fibration given in Theorem 0.2.

Proof: Assume $f:X\rightarrow S$ is given as in Lemma 1.2
and $j:J\rightarrow S$
as given in Lemma 1.6.
We have $K_J=j^*(K_S+\Delta)$, where $\Delta=\Delta_{X/S}=\Delta_{J/S}$, by
Lemma 1.6 and the fact that $\Pic S$ has no torsion as $S$ is rational.
Now (continuing the notation of the proof of
Theorem 0.2)
$$K_S+\Delta=\sum (c_i-(m_i-1)/m_i)E_i=\sum (d_i+a_i+b_i')E_i,$$
and by (2.2) and (2.3), $d_i+a_i+b_i'$ is a non-negative integer.
Hence $K_S+\Delta$ is an effective divisor, and so is $K_J$. Thus
$K_S+\Delta$ has a Zariski decomposition, which by uniqueness is
$$K_S+\Delta=0+\sum (d_i+a_i+b_i')E_i.$$
Thus to obtain a minimal model for $J$, we contract those curves $E_i$
on $S$ for which $d_i+a_i+b_i'>0$, to obtain a surface $S'$, and take
$j:J'\rightarrow S'$ to be a relatively minimal model of $J\rightarrow
S'$. Since $S\rightarrow S'$ contracts a subset of the exceptional divisors
of $\phi:S\rightarrow \bar S$, there exists a morphism $S'\rightarrow \bar S$.

To show $J'$ is Calabi-Yau, we now only need to show that $\chi(\O_J)=0=
h^1(\O_J)=h^2(\O_J)$. But this follows from the corresponding facts
for $X$ by [4, 2.3, 2.4].
$\bullet$

\proclaim Proposition 2.3. If $f:X\rightarrow S$ is an elliptic fibration
with $X$ birational to a Calabi-Yau threefold, and $\bar f:\bar X\rightarrow
\bar S$ the minimal model constructed in Theorem 0.2,
then $\bar S$ is either rational
or is birational to an Enriques surface. In the latter case, the minimal
resolution of $\bar S$ is a minimal Enriques surfaces.

Proof: It is clear that $\kappa(S)\le 0$, and if $\kappa(S)=-\infty$, then $S$
must be rational since $h^1(\O_X)=0$. If $S$ has Kodaira dimension zero,
then it must be a K3 or Enriques surface for the same reason. If it is a
K3, then by [4], Prop. 2.2, $h^2(\O_X)\ge 1$ since $h^2(\O_S)=1$, and
this is again a contradiction. Thus $S$ must be an Enriques surface.
By [4], Theorem 3.1 b), the minimal resolution of $\bar S$ is a minimal
Enriques surface. $\bullet$

\proclaim Definition 2.4.
Let $S$ be a normal Gorenstein surface. An irreducible
curve $C$ in $S$ is an {\it exceptional curve of the first kind} if
$K_S.C<0$ and $C^2<0$. $S$ is said to be
minimal if $S$ contains no exceptional curves
of the first kind.

\proclaim Theorem 2.5. Let $S$ be a normal Gorenstein surface.
There exists a morphism $S\rightarrow S'$ with
$S'$ minimal. If $K_{S'}$ is not nef, then either (i) $-K_{S'}$ is
numerically ample (i.e. $-K_{S'}.C>0$ for all irreducible curves
$C$ on $S'$) and $\rho(S')=\dim\Pic(S')\otimes{\bf Q}=1$
or (ii) there is a smooth curve $C$ and a
map $S'\rightarrow C$ such that no fibre contains an exceptional curve
of the first kind.  In this case all fibres of $\pi$ are irreducible.
In any event, $S'$ is called a minimal model of $S$.

Proof: [22, Thm. 4.9 and Lemma 4.6] The theorem is proved by contracting
exceptional curves of the first kind until there are none left,
and then analyzing the case that $K_{S'}$ is not nef.
$\bullet$

{\it Remark 2.6.} If $S$ has only DuVal singularities, then an
exceptional curve $C$ of the first kind is either a $-1$-curve disjoint
from the singular locus, or else passes through precisely one singularity
of type $A_n$. Furthermore, if $\tilde S\rightarrow S$ is the minimal
resolution of $S$, we can write the resolution of the $A_n$ singularity
as a union of irreducible curves $C_1,\ldots,C_n$, $C_i.C_{i+1}=1$, $
1\le i \le n-1$. Then the proper transform of $C$ on $\tilde S$ is a $-1$-curve
and
intersects
only $C_1$ or $C_n$. In particular, $C$ then contracts to a smooth point.
[22, Example 1.2.]

\proclaim Definition 2.7. A Gorenstein log Del Pezzo surface is a surface
$S$ with only DuVal singularities with $-K_S$ ample. The rank of $S$
is $\rho(S)=\dim \Pic(S)\otimes{\bf Q}$.

Clearly if $\rho(S)=1$ then $-K_S$ is ample if it is numerically ample.

\proclaim Proposition 2.8. Let $S$ be a rational surface with only $A_n$
singularities. Then there exists a minimal model $S'$ of $S$ which is
either a (i) rank one Gorenstein
log Del Pezzo surface, or (ii) there is a map $S'\rightarrow
\Pone$. Each fibre is a $\Pone$, and each reducible fibre of
$\tilde S'\rightarrow\Pone$ where $\tilde S'$ is a minimal
desingularization of $S'$ either has three components $C_1,C_2,C_3$ with
intersection matrix
$$\bordermatrix{&C_1&C_2&C_3\cr
C_1&-2&1&0\cr
C_2&1&-1&1\cr
C_3&0&1&-2\cr}$$
or four components $C_1,C_2,C_3,C_4$ with
$$\bordermatrix{&C_1&C_2&C_3&C_4\cr
C_1&-2&1&0&0\cr
C_2&1&-2&1&1\cr
C_3&0&1&-2&0\cr
C_4&0&1&0&-1\cr}.$$
Here in the first case $C_1$ and $C_3$ are contracted to obtain $S'$, giving
two $A_1$ singularities, and in the second case $C_1,C_2,$ and $C_3$
are contracted to give one $A_3$ singularity. We call $S'$ a minimal
ruled surface, and the singular fibres of this surface are the fibres
containing singularities of the surface.

Proof:  Let $\pi:\tilde S'\rightarrow S'$ be a minimal
resolution of singularities.
By Remark 2.6, $S'$ has only $A_n$ singularities, so $K_{\tilde S'}=
\pi^*K_{S'}$, and since $\tilde S'$ is rational, $K_{\tilde S'}$ is not nef,
and hence $K_{S'}$ is not nef. By Theorem 2.5, we then have the two cases
given.

All that remains to check is the possible reducible fibres of the fibration
$\tilde S'\rightarrow \Pone$ if $S'$ is a minimal ruled surface. Let
$C_1,\ldots,
C_r$ be the irreducible components of a given fibre of this fibration, with
$C_r$ the one curve  of the fibre not contracted by $\pi$. Each $C_i$ is
a $\Pone$, and $C_i^2=-2$, $i<r$, and $C_r^2=-1$, since the singularities
are DuVal, the resolution is minimal, and there must be at least one
component with self-intersection $-1$. $C_r$ cannot intersect more than
two other components, as a contraction of $C_r$ would then yield three
components meeting at a smooth point, and the fibre must always be a tree
of $\Pone$'s. If $C_r$ meets two other curves, then after contracting $C_r$
we obtain two curves which intersect, each with self-intersection
$-1$. Contracting one of these yields
a curve of self-intersection zero. By Zariski's Lemma, this must then be
the only component, $r=3$, and we are in the first case.

If $C_r$ meets only one curve, say $C_{r-1}$, then we can contract $C_r$,
and now $C_{r-1}$ is a $-1$ curve. We repeat the process, and for this
to terminate, we must eventually arrive at the previous situation. Thus the
singularity would be $D_{r-1}$ unless $r=4$, in which case the singularity
is $A_3$, and we are in the second case. $\bullet$

\proclaim Theorem 2.9. Let $f:X\rightarrow S$ be an elliptic fibration with
$X$ a Calabi-Yau threefold and $S$ rational, and let $j:J\rightarrow S$
be its jacobian fibration, and $\bar S$, $S'$ as in Proposition 2.2.
Then
$\bar S$ has only $A_n$ singularities, and a minimal model for $\bar S$,
$\bar S'$, is
either $\Ptwo$, $F_e$, $0\le e \le 12$, a rank one Gorenstein
log Del Pezzo surface with
only $A_n$ singularities, or else is a minimal ruled surface with only
$A_n$ singularities and at most
four singular fibres. There is a birational
morphism $S'\rightarrow \tilde S'$, where
$\tilde S'$ is the minimal resolution of $\bar S'$, and there is a
birational
morphism $\tilde S'\rightarrow F_2$ if $\bar S'$ is a log Del Pezzo surface,
and
a birational
morphism $\tilde S'\rightarrow F_e$, $e=0,1$ or $2$ if $\tilde S'$ is
a minimal ruled surface.
$$\matrix{S&\mapright{\pi}&S'&\mapright{\pi'}&\tilde S'\cr
&&\mapdown{}&&\mapdown{}
\cr
&&\bar S&\mapright{}&\bar S'\cr}\leqno{(2.4)}$$

Proof: That $\bar S$ has only $A_n$ singularities is Proposition 3.3,
and hence $\bar S'$ has only $A_n$ singularities by Remark 2.6. First
suppose that $\bar S'$ is non-singular. Then there is a birational
morphism $S'\rightarrow
\bar S'$ by Proposition 2.2.
By [4], Prop. 3.5, $\bar S'$ is then either $\Ptwo$ or $F_e$,
$0\le e\le 12$.

If $\bar S'$ is singular, then
since there is a morphism $S'\rightarrow \bar S'$
by Proposition 2.2, this map must factor through the minimal resolution
$S'\rightarrow \tilde S' \rightarrow \bar S'$.

If $\bar S'$ is a log Del Pezzo surface, then by [17], Lemma 3, there is a
morphism $\tilde S'\rightarrow F_2$.

If $\bar S'$ is a minimal ruled surface, consider the diagram
$$\matrix{S&\mapright{\pi}&S'&\mapright{\pi'}&\tilde S'&\mapright{\pi''}&\tilde
S''\cr
&&\mapdown{}&&\mapdown{}&&\mapdown{}\cr
&&\bar S&\mapright{}&\bar S'&\mapright{}&\Pone\cr}$$
where $\pi''$ is the contraction of those $-1$-curves in reducible fibres
of $\tilde S'\rightarrow\Pone$
which have
four components (i.e. $C_4$ in the notation of Proposition 2.8), so that
each reducible fibre of $\tilde S''\rightarrow\Pone$ has
precisely three components. (This
contraction is merely to avoid having two seperate cases to analyze.)

Now the discriminant locus of $j:J\rightarrow S$, $\Sigma\sim 12\Delta$,
can be written
as
$$\Sigma=\Sigma_0+J_{\infty},$$
where $J_{\infty}$ is the fibre of the $J$-morphism to $\Pone$ over $\infty$,
and $\Sigma_0=\sum 12a_iD_i$ as in Lemma 1.2, 2). By Cor. 2.1, we can
find a point $P\in \Pone$ such that $J_P$ is disjoint from any curve on
$S$ contracting to a singular point on $\bar S$, and then
$$\Sigma\sim \Sigma_0+J_P.$$
Replace $\Sigma$ by the linearly equivalent divisor $\Sigma_0+J_P$, so that
$\Sigma\subseteq S$ is disjoint from all curves contracted to singular points
on $\bar S$.
Now $\pi_*\Sigma=\Sigma'\sim -12K_{S'}$ (as $\pi_*\Delta=-K_{S'}$)
and so $\pi''_*\pi'_*\pi_*\Sigma=\tilde\Sigma''\sim -12K_{\tilde S''}$,
and $\tilde\Sigma''$ is disjoint from the curves $C_1$ and $C_3$ on any
reducible fibre of $\tilde S''\rightarrow\Pone$.

If $\sigma\subseteq \tilde S''$ is a section of $\tilde S''\rightarrow \Pone$,
with $\sigma^2\le -3$, then $\sigma.K_{\tilde S''}\ge 1$,
$\sigma.\tilde\Sigma''
<0$, and so
$\sigma\subseteq\tilde\Sigma''$.
 But $\sigma$ must pass through either the $C_1$ or the $C_3$
component of a singular fibre of $\tilde S''\rightarrow \Pone$, which is
then a contradiction.
Thus $\tilde S''$ does not
contain a section of self-intersection $\le -3$, so there is a morphism
$\tilde S''\rightarrow F_e$, $e=0,1$ or $2$.

Now let
$C$ be an irreducible component of $\tilde\Sigma''$ which
dominates $\Pone$, and let
$\tilde C$ be its normalization. As above, $C$ cannot be a section.
We can take a basis of $\Pic \tilde S''$
to be $\sigma_0$, a section of $\tilde S''\rightarrow\Pone$ with $\sigma_0^2
=-e$, $f$ the class of a fibre and $C_2^i,C_3^i$, $1\le i\le n$,
where $C_2^i$ and $C_3^i$ are the components of the $n$ reducible fibres, using
the numbering in Proposition 2.8.
We can assume that $\sigma_0.C_3^i=0$ for all $i$ by interchanging
$C_1^i$ and $C_3^i$ if necessary.
We can then write
$$C\sim a\sigma_0+bf+\sum_{1\le i \le n,j=2,3} c^i_j C^i_j,$$
with $a\ge 2$, since $C$ is not a section.
Now as remarked above, $C.C_1^i=C.C_3^i=0$, for all $i$, which
yields the equalities
$$\eqalign{a+c_2^i&=0\cr
c_2^i-2c_3^i&=0.\cr}$$
Hence
$$C\sim a\sigma_0+bf-\sum_{1\le i \le n} (aC_2^i+{a\over 2}C_3^i).$$
We also have
$$K_{\tilde S'}\sim -2\sigma_0-(2+e)f+\sum_{1\le i \le n} (2C_2^i+C_3^i).$$
By adjunction, we then obtain
$$2p_a(\tilde C)-2\le -a^2e+ae+2ab-2a-2b+n(-{a^2\over 2}+a),$$
and Riemann-Hurwitz tells us that
$$2p_a(\tilde C)-2=-2a+R,$$
where $R$ is the degree of the ramification divisor of
$\tilde C\rightarrow\Pone$. It is easy to see that $R\ge an/2$, since
if $P_1,\ldots,P_n$ are the points of $\Pone$ where the fibres of
$\tilde S''\rightarrow \Pone$ are reducible, then every point of $\tilde C$
over $P_i$ is a ramification point of $\tilde C\rightarrow \Pone$.
Thus we obtain the inequality
$$-a^2e+ae+2ab-2a-2b+n(-{a^2\over 2}+a)\ge -2a +an/2,$$
which reduces to (remembering that $a\ge 2$)
$${b \over a}\ge {n\over 4}+{e\over 2}.$$
Thus if $n>4$, we see that ${b\over a}>1+e/2$.
But
$$\tilde\Sigma''\cong -12K_{\tilde S'} \cong a'\sigma_0+b'f+\sum
c_i'C_2^i+d_i'C_3^i,$$
with $b'/a'=1+e/2$, so no sum of components with $b/a>1+e/2$ can yield
such a $\tilde\Sigma''$. Thus $n\le 4$.
$\bullet$

\proclaim Proposition 2.10. Let $f:X\rightarrow S$ be an elliptic fibration
with $X$ Calabi-Yau and $S$ birational to an Enriques surface, and let
$\bar f:\bar X\rightarrow \bar S$ be as in Theorem 0.2. Then there exists
a surface $T$ which is either a K3 or abelian surface, and a map
$g:T\rightarrow
\bar S$, where $g$ is a Galois covering which is \'etale in codimension one,
with $\bar X\times_{\bar S} T$ birational to $E\times T$ over $T$ for some
elliptic curve $E$.

Proof: This follows immediately from [19], Appendix, Thm. (2.1). $\bullet$

{\hd \S 3. Tate-Shafarevich Groups.}

We first review the notion of the Tate-Shafarevich group for elliptic
threefolds introduced in [3]. See
[3] \S\S 1 and 2 for more details.

Let $A$ be an elliptic curve with a rational point defined
over a field $K=K(S)$, the function field of
the surface $S$. Let $i:\eta\rightarrow S$ be the inclusion of the
generic point.
We define the Weil-Ch\^atelet group
$$WC(A):=H^1(\eta,A),$$
thinking of $A$ as a sheaf in the \'etale topology. This is the group of
torsors over $A$, i.e. the set of curves of genus 1 over $K$ with jacobian
$A$. Now for each $s\in S$, there is a natural localization map
$$loc_{\bar s}:WC(A)\rightarrow WC(A_{\bar s})$$
where $A_{\bar s}=A\times_{\eta} \eta_{\bar s}$ and $\eta_{\bar s}$ is
the function field of $\O_{S,\bar s}$,
the strict henselization of the local ring of $S$
at $s$. The map is given by $E\mapsto E\times_{\eta} \eta_{\bar s}$.
We define the Tate-Shafarevich group to be
$$\TS_S(A):=\bigcap_{s\in S}\ker(loc_{\bar s}).$$
This turns out to be also $H^1(S,i_*A)$, with cohomology in the \'etale
topology.

We can identify $\TS_S(A)$  with the set
$$\{E\in WC(A)|\hbox{$X_E\rightarrow S$ has a rational section locally
in the \'etale topology at every point $s\in S$}\}$$
where $X_E\rightarrow S$ is proper and a model of $E$ (i.e. $(X_E)_{\eta}
=E$.) Thus if $E\in\TS_S(A)$, $f:X\rightarrow S$ a model for $E$, then $f$ can
have only isolated multiple fibres. In this case, if
$s \in S$ and $X_s$ is an isolated
multiple fibre, $f$ still has a rational section locally at $s$. Such
an isolated multiple fibre is called a {\it locally trivial isolated
multiple fibre.} (See [3], 2.19.)

If $f:X\rightarrow S$ is an elliptic fibration with only locally trivial
isolated multiple fibres, with $X$ birational to a Calabi-Yau, then we have
seen
in Proposition 2.2 that its jacobian, $j:J\rightarrow S$, is also
birational to a Calabi-Yau. Thus as a first approximation to a finiteness
theorem, we would hope that $\TS_S(A)$ is finite where $A$ is the generic
fibre of $j$.

In [3], \S 1, we have shown how to calculate $\TS_S(A)$ if $j:J\rightarrow
S$ is a Miranda model.
{}From this we obtain almost immediately the following crucial theorem. Note
that
this fails for elliptic K3 surfaces, since $h^2(\O_X)\not=0$ when $X$
is K3. This explains why our main theorem is not true in the K3 case.

\proclaim Proposition 3.1. If $J\rightarrow S$ is a Miranda fibration with
generic fibre $A$ and
$J$ a Calabi-Yau threefold and $S$ projective,
then $\TS_S(A)$ is a finite
group.

Proof: By [3], Theorem 2.24, $\TS_S(A)$ is an extension of $(\qz)^r$ by
a finite group, where $r$ is the corank of $\TS_S(A)$, and
$$r=b_2(J)-\rho(J)-(b_2(S)-\rho(S))$$
where $\rho$ is the rank of the Picard group and $b_2$ is the second
Betti number. Since $h^2(\O_J)=h^2(\O_S)=0$,
$b_2(J)=\rho(J)$ and $b_2(S)=\rho(S)$. Thus $r=0$ and $\TS_S(A)$ is a finite
group. $\bullet$

Of course, an elliptic Calabi-Yau threefold might have multiple fibres which
are not locally trivial. The group $\TS_{S-Z}(A)/\TS_S(A)$ measures the
additional fibrations one obtains by allowing multiple fibres along $Z$.
Since only isolated multiple fibres are allowed on the minimal model
of an elliptic Calabi-Yau threefold, we consider sets $Z$ which
can be contracted to
a set of points. This motivates the following proposition. The proof
requires familiarity with the results and notation of [7].

\proclaim Proposition 3.2. Let $j:J\rightarrow S$ be a Miranda model with
generic fibre $A$,
and let $E_i\subseteq S$, $1\le i \le n$ be irreducible projective
curves with the intersection matrix $(E_i.E_j)$ negative definite. Then the
group
$$\TS_{S-\bigcup E_i}(A)/\TS_S(A)$$
is a finite group.

Proof: Put $Z=\bigcup E_i$, and order the $E_i$'s so that
$E_1,\ldots,E_{n_1}$ are of fibre-type $I_0$  with $j^{-1}(E_i)=E_i\times C_i$
for elliptic curves $C_i$, $1\le i \le n_1$,
$E_{n_1+1},\ldots,E_{n_2}$
are of fibre type $I_a$, Case $I$,
$a\ge 1$.
(See [7], \S 1, for the distinction between Case $I$
and Case $I^*$ for curves of fibre type $I_a$, $a\ge 1$.)
The remaining components of $Z$ will be curves not
of this type. Then from [7, 2.11],
$$H^2_Z(S,i_*A)\subseteq \bigoplus_{i=1}^{n_1} MW(j^{-1}(E_i)/E_i)_{tors}
\oplus\bigoplus_{i=n_1+1}^{n_2}\qz\oplus G,$$
where $G$ is a finite group and $MW(j^{-1}(E_i)/E_i)$ is the Mordell-Weil
group of the elliptic surface $j^{-1}(E_i)\rightarrow E_i$.
Indeed, from [7], 2.11, $G$ has various finite contributions
from curves of type $I_a$, Case $I^*$, a finite number of
collisions of type $I_0^*+III$, and the torsion parts
of Mordell-Weil groups of non-trivial elliptic surfaces, which by the
Mordell-Weil theorem
are finitely generated, and hence have finite torsion.

By the exact sequence of local cohomology, we have
$$\TS_{S-Z}(A)/\TS_S(A)\cong \ker(H^2_Z(S,i_*A)\mapright{\phi} H^2(S,i_*A)).$$
So, while $H^2_Z(S,i_*A)$ may well be infinite, we can show that the kernel
of $\phi$ is finite. Recall from [7] that if $\gamma\in H^2_Z(S,i_*A)$,
then $\phi(\gamma)=0$ if and only if the so-called first and second
obstructions
of
$\gamma$ are zero. See [7], \S 4, for details.

Let $\gamma_1,\ldots,\gamma_{n_2}$ be the invariants of an element $\gamma$ of
$H^2_Z(S,i_*A)$ along the curves $E_1,\ldots, E_{n_2}$, with $\gamma_1,
\ldots,\gamma_{n_1}
\in (\qz)^{\oplus 2}=MW(j^{-1}(E_i)/E_i)_{tors}$ and
$\gamma_{n_1+1},\ldots,\gamma_{n_2}\in \qz$.
(See [7], Remark 2.12)
Consider the first obstruction along the curves $E_{n_1+1},\ldots,E_{n_2}$.
By [7] Theorems 4.5 and 4.3,
this obstruction being zero requires in particular that for all $l$
$$\left(\sum_{i=n_1+1}^{n_2} \gamma_iE_i.E_l\right)+k_l=0,\leqno(3.1)$$
where $k_l\in\qz$ depends on the other invariants of
$\gamma$ along the curves $E_i$, $i>n_2$. These invariants take
values in the finite
group $G$, and thus $k_l$ can only take on a finite number of possible values.
Now the intersection matrix $(E_i.E_l)_{n_1+1\le i,l \le n_2}$ is negative
definite, and as a result, the number of solutions to (3.1) is finite for any
given set of values for the $k_l$'s. Thus there are only a finite number
of possibilities for the $\gamma_i$, $n_1+1\le i \le n_2$.

A similar argument works for the $\gamma_i$, $1\le i \le n_1$, but we need to
use
the second obstruction, which is a codimension 2 algebraic cycle on
$J$. Identifying
$(\Pic C_i)_{tors}$ with
$(\qz)^{\oplus 2}$, following [7], Examples 5.2, 1) and 2),
we obtain a cycle
of the following form. If $\pi_i:E_i\times C_i\rightarrow C_i$ is the
second projection, $f_i:E_i\times C_i\rightarrow J$ the inclusion,
$c_i\in (C_i)_{tors}$, then the second obstruction
is a cycle
$$\alpha=\left(\sum_{i=1}^n f_{i*}\pi_i^*(\gamma_i)\right)+\beta,$$
where $\beta$ is a cycle which again only takes on a finite number of possible
different values. In order to realise this set of invariants, we need
$\alpha$ to be rationally equivalent to zero. In order for this to be possible,
in particular we must have $\pi_{l*}f_l^*(\alpha)$ linearly equivalent to zero
on each curve $C_l$, i.e.
$$\pi_{l*}f_l^*(\alpha)=\left(
\sum_{i=1}^n \gamma_iE_i.E_l\right)+k_l=0,\leqno{(3.2)}$$
as an equation in $(\qz)^{\oplus 2}$,
where the $k_l$ come from the cycle $\beta$ and thus take on only a finite
number of possible values. We then argue as before, and the $\gamma_i$,
$1\le i \le n$, take on only a finite number of values. $\bullet$

As an application of the above proof, we obtain

\proclaim Proposition 3.3. The singularities occuring in $\bar S$ in Theorem
0.2 are all $A_n$ singularities, i.e the minimal resolution consists of a chain
of $-2$ curves.

Proof: As usual, we have $f:X\rightarrow S$, $j:J\rightarrow S$ as given in
Lemmas 1.2 and 1.6.
We have the resolution $S\rightarrow \bar S$ of the singularities
of $\bar S$. Let $P\in\bar S$ be a singular point. Then
$\phi^{-1}(P)$ consists of curves $E_1,\ldots,E_r$ which are all isomorphic
to $\Pone$. By Corollary 2.1, there are two cases. Either $j^{-1}(E_i)
=E_i\times C_i$ for an elliptic curve $C_i$, or all $E_i$ are of fibre
type $I_a$, $a\ge 1$ and Case $I$. Thus either formulae (3.1) or (3.2)
are relevant, except that there is no $k_l$. Thus we obtain,
with $\gamma_1,\ldots, \gamma_r\in (\qz)^{\oplus 2}$ or $\qz$ depending
on the two cases, the following system of equations in $(\qz)^{\oplus 2}$
or $\qz$:
$$\sum_{i=1}^r \gamma_iE_i.E_l=0,\quad 1\le l \le r.$$
The multiplicity of the fibres over $E_i$ in the fibration $f:X\rightarrow S$,
$m_i$, is the order of $\gamma_i$ ([7], 2.12).

Now in the notation of the proof of Theorem 0.2, $a_i=b_i'=0$ for all $i$. Thus
by (2.3), if $d_i=0$ for some $i$, then we must have $m_i>0$, and hence
$\gamma_i\not=0$. Let $\C=(E_i.E_j)$ be the intersection matrix
of the $E_i$'s. Then $\gamma_i=0$ if the $i$th row of $\C^{-1}$ consists
of integers. If the resolution is minimal, it is easy to check that there
exists integral rows in $\C^{-1}$ for singularities $D_m$ or $E_m$. (In fact,
for $E_8$, $\C^{-1}$ itself is integral). If the resolution is not minimal,
then it is a blowup of a minimal resolution, and it is easy to see that
blowing up does not effect the integrality of a row in $\C^{-1}$. $\bullet$

{\it Example 3.4.} Let $f_1:Y_1\rightarrow \Pone$ and $f_2:Y_2\rightarrow
\Pone$
be two rational elliptic surfaces, with $Y_1\rightarrow \Pone$ having two
singular fibres of type $I_0^*$. Such a surface exists by
[16], Theorem 4.1. Let $X=Y_1\times_{\Pone}
Y_2$ with $p_1,p_2:X\rightarrow Y_1,Y_2$ the projections. If $Y_2$
is chosen so that no singular fibres of $f_1$ and $f_2$ map to the same point
in $\Pone$, then $X$ is a non-singular Calabi-Yau threefold, and $p_1$ and
$p_2$ are elliptic fibrations. (See [23], \S 2)
Let $P\in\Pone$ be a point with the fibre $(Y_1)_P$
of type $I_0^*$ and $(Y_2)_P=E\subseteq Y_2$ a non-singular elliptic
curve. $(Y_1)_P$ has five components $C_1,\ldots,C_5$ with $C_5$ having
multiplicity 2. We have $p_1^{-1}(C_i)\cong C_i\times E$. We wish to
perform a logarithmic transformation along $C_1,\ldots, C_4$.

Let $P_1,P_0\in E\subseteq Y_2$ be two points with $P_1-P_0$ a 2-torsion
element in $\Pic E$, and consider the cycle on $X$ given by $p_2^*(P_1-P_0)$.
Since $Y_2$ is rational, $P_1-P_0$ is rationally equivalent to zero on
$Y_2$, so $p_2^*(P_1-P_0)\sim 0$ on $X$.
Now
$$\eqalign{p_2^*(P_1-P_0)&=
\sum_{i=1}^4 \alpha_{i*}{(P_1\times C_i-
P_0\times C_i)}+2\alpha_{5*}(P_1\times C_5-P_0\times C_5)\cr
&\sim
\sum_{i=1}^4 \alpha_{i*}{(P_1\times C_i-
P_0\times C_i)}\cr}$$
where $\alpha_i:C_i\times E\rightarrow X$ is the inclusion. (The last
rational equivalence is because $P_1-P_0$ is 2-torsion in $E$.)
Because this expression is rationally equivalent to zero,
as in [7], Example 5.2. we can
then obtain a threefold $X'\rightarrow Y_1$ whose jacobian is $p_1:X\rightarrow
Y_1$ with fibres of multiplicity 2
along $C_1,\ldots,C_4$. $X'$ is not minimal, but as in the proof of Theorem
0.2,
one can contract $C_1,\ldots,C_4$ to obtain a surface $\bar Y_1$ and a minimal
Calabi-Yau threefold $\bar X'$ with a fibration $\bar X'\rightarrow \bar
Y_1$. $\bar Y_1$ has four $A_1$ points.

We now put results of the last section and this section together. This
statement is very close to our final result; we will only need to show
how Ogg-Shafarevich theory works in families, which we will do in the
next section.

\proclaim Proposition 3.5. Let $f:X\rightarrow S$ be a Calabi-Yau
elliptic fibration with $S$ rational and
generic fibre $E$ and $Jac(E)=A$. Let $\bar S'$ and
$\tilde S'$ be as in Theorem 2.9. Then there exists sections
$a\in\Gamma(\tilde S',\omega_{\tilde S'}^{-4})$ and
$b\in\Gamma(\tilde S',\omega_{\tilde S'}^{-6})$ such that the Weierstrass
model
$$W(\omega_{\tilde S'}^{-1},a,b)\rightarrow \tilde S'$$
is minimal and has generic fibre $A$.
If $\Sigma_{red}$ is the reduced discriminant locus of this Weierstrass
model, let $Z\subseteq \tilde S'$ be the union of $sing(\Sigma_{red})$ and
the exceptional locus of $\tilde S'\rightarrow \bar S'$. Then
$E\in\TS_{\tilde S'-Z}(A)$. Furthermore, this latter group is finite.

Proof. As in Proposition 2.2, we have a jacobian fibration $j':J'\rightarrow
S'$ with $K_{J'}=0$. We obtained a composed fibration $J'\rightarrow \tilde
S'$.
This fibration has a rational section, and hence we can find a birationally
equivalent minimal Weierstrass fibration
$$w:W:=W(\L,a,b)\rightarrow \tilde S'.$$
It is then easy to see that $\L\cong \omega_{\tilde S'}^{-1}$. Indeed,
$K_W=w^*(\omega_{\tilde S'}\otimes\L),$ and if $T=
\tilde S'-sing(\Sigma_{red})$, $W_T=W\times_{\tilde S'} T$
has a canonical crepant
resolution of singularities $\phi:\tilde W_T\rightarrow W_T$ with
$K_{\tilde W_T}=\phi^*K_{W_T}$. Furthermore, $\tilde W_T\rightarrow T$
and $J'_T=J'\times_{\tilde S'}T\rightarrow T$ are birationally equivalent
and both are relatively minimal; since $K_{J'_T}=0$, this implies that
$K_{\tilde W_T}=0$ also, and hence $K_{W_T}=K_W=0$. Thus $\L\cong\omega_{\tilde
S'}^{-1}$.

Now consider the fibration $\bar X\rightarrow \bar S'$. $\Sigma^m$ for this
fibration consists of a finite number of points, as $\Sigma^m$ for $\bar X
\rightarrow \bar S$ consists of isolated points. Every one of these points
must be either at a singular point of $\bar S'$ or a singular point of the
reduced discriminant locus of
$\bar X\rightarrow \bar S$, $\bar\Sigma_{red}$, by [3], Cor. 3.2. Let $Z'
=sing(\bar\Sigma_{red})\cup sing(\bar S')$. Then
$\bar S'-Z'=\tilde S'-Z$, and clearly
$E\in \TS_{\bar S'-Z'}(A)=\TS_{\tilde S'-Z}(A)$.

Finally, to show the finiteness statement, let $J_M\rightarrow M$
be a Miranda model for $A$ obtained by blowing up $\tilde S'$. The
exceptional locus of $\psi:M\rightarrow \tilde S'$ maps to $Z$, since
we only need to blow-up points of $sing(\Sigma_{red})$. Now $\psi^{-1}(Z)$
consists of a finite number of points and a union $E_1\cup\cdots\cup E_n$ of
$\Pone$'s whose intersection matrix is negative definite by Grauert's
criterion. Now Proposition 3.1 tells us that $\TS_M(A)$ is finite, Proposition
3.2 tells us that $\TS_{M-E_1\cup\cdots\cup E_n}(A)$ is then finite,
and [3], Cor. 3.2 tells us that removing a finite number of further
points can only increase the Tate-Shafarevich group by a finite group,
so that $\TS_{M-\psi^{-1}(Z)}(A)=\TS_{\tilde S'-Z}(A)$ is a finite group.

{\hd \S 4. Proof of the Main Theorem.}

\proclaim Definition 4.1. A family of elliptic three-folds is a triple
$(\X,\SS,\T)$ of non-singular quasi-projective varieties with a diagram
$$\matrix{\X&&\cr
\mapdown{f}&\mapse{\pi}&\cr
\SS&\mapright{g}&\T\cr}$$
with $f,g$ and $\pi$ projective, $g$ and $\pi$ smooth of relative dimension
2 and 3 respectively, and the generic fibre of $f$ a geometrically regular
curve of genus 1.

\proclaim Lemma 4.2. Let $g:\SS\rightarrow \T$ be a family of
smooth projective surfaces, $\L$ a line
bundle on $\SS$. Then there exists a finite number of families
$(\X_i,\SS_i,\T_i)$ of elliptic fibrations with $f_i: \X_i\rightarrow
\SS_i$ having a section, along with $\X_i\rightarrow
W_i(\L_i,a_i,b_i)$ a resolution of singularities of a Weierstrass
model over $\SS_i$, and also with Cartesian
diagrams
$$\matrix{\SS_i&\mapright{}&\T_i\cr
\mapdown{\alpha_i}&&\mapdown{}\cr
\SS&\mapright{}&\T\cr}$$
with $\L_i\cong \alpha_i^*\L$, such that the following property holds. For any
$t\in \T,$ $S=\SS_t$,
$a\in\Gamma(S,\L^{\otimes 4}|_S)$,
$b\in\Gamma(S,\L^{\otimes 6}|_S)$ with $4a^3+27b^2\not=0$, there exists
a $t'\in \T_i$ for some $i$ with $a_i|_{(\SS_i)_{t'}}=a$, $b_i|_{(\SS_i)_{t'}}
=b$.
Furthermore, we can assume that
$f_i:\X_i\rightarrow \SS_i$ is flat away from the singularities of the
reduced discriminant locus of $f_i$.

Proof: Let $\E=\L^{\otimes 4}\oplus\L^{\otimes 6}$, and let
$\T'={\bf V}((g_*\E)^{\vee})$ and $\SS'={\bf V}((g^*g_*\E)^{\vee})$, where
${\bf V}(\F)={\bf Spec}(S(\F))$, $S(\F)$ the symmetric algebra of
$\F$. By semicontinuity, we can split $\T$ up into a finite number of locally
closed subsets of $\T$ on which $g_*\E$ is locally free, so we can
assume that $g_*\E$ is indeed locally free. We have a diagram
$$\matrix{\SS'&\mapright{\pi'}&\SS\cr
\mapdown{g'}&&\mapdown{g}\cr
\T'&\mapright{\pi}&\T\cr}$$
Now $\pi'^*\E$ has a universal section given by the composition
$$\O_{\SS'}\rightarrow \pi'^*g^*g_*\E\rightarrow \pi'^*\E$$
where the first map is the universal section of $\pi'^*g^*g_*\E$
([9], 9.4.9). If we write this section as $(a,b)$, $a\in\Gamma(\SS',
\pi'^*\L^{\otimes 4})$ and $b\in\Gamma(\SS',\pi'^*\L^{\otimes 6})$, then if
$t'\in\T'$ corresponds to a point $(t,a',b')$, $t\in\T$, $a'\in\Gamma(\SS_t,
\L^{\otimes 4}|_{\SS_t})$, $b'\in\Gamma(\SS_t,\L^{\otimes 6}|_{\SS_t})$, then
$a,b|_{(\SS')_{t'}}=a',b'$.
We then have a universal Weierstrass model
$$W'=W(\pi'^*\L,a,b)\rightarrow \SS',$$
and we can replace $\T'$ with the open subset
$$\{t'\in\T'|\hbox{$4a^3+27b^2$ is not identically zero on $(\SS')_{t'}$}\},$$
and replace $\SS'$ with the inverse image of this open subset. Thus we can
assume that no fibre of $g'$ is contained in the discriminant locus
 $\Sigma\subseteq
\SS'$.

Now let $\X'\rightarrow W'$ be a resolution of singularities, and let
$\T_1\subseteq \T'$ be a dense open set on which $\X'\rightarrow\T'$ is smooth,
$\SS_1=\SS'\times_{\T'}\T_1$,
$\X_1=\X'\times_{\T'}\T_1$,
and
$$W_1=W(\pi'^*\L|_{\SS_1},a|_{\SS_1},b|_{\SS_1})\rightarrow \SS_1.$$
As in [15], \S 7, this resolution can be performed
canonically away from the singular points of the reduced discriminant locus
$\Sigma_{red}$ of $W'\rightarrow \SS'$. Thus we can assume that
$\X'\rightarrow\SS'$ is flat away from $sing(\Sigma_{red})$.
Note also that a Weierstrass model is always non-singular at the section
at infinity, so the resolution can be performed so that $\X'\rightarrow
\SS'$ has a section.
Now replace $\T'$ with a finite disjoint union of non-singular locally
closed subsets of $\T'-\T_1$, whose union is $\T'-\T_1$, replace $W'$ and
$\SS'$ with restrictions of these to the new $\T'$, and then construct
a resolution of singularities of the smaller $W'$. By Noetherian induction,
we eventually cover all of $\T'$ in this fashion, proving the theorem.
$\bullet$

\proclaim Theorem 4.3. Let $(\bar\X,\bar\SS,\T)$ be a family of
elliptic fibrations such that $\bar f:\bar\X\rightarrow \bar\SS$ has
a section, and let $\SS\subseteq\bar\SS$ be an open subset,
$\X=\SS\times_{\bar\SS}\bar\X$. Assume furthermore that the
reduced
discriminant locus $\Sigma_{red}$
of $f:\X\rightarrow\SS$ is smooth over $\T$, every component of
$\Sigma_{red}$ surjects onto $\T$, and that
$f$ is flat.
Furthermore, assume that
for all $t \in\T$,
$\TS_{\SS_t}(A_t)$ is a finite group,
where $A_t$ is the generic fibre of
$\X_t\rightarrow \SS_t$. Then there exists a finite collection of families
of elliptic threefolds $(\bar\X_i,\bar\SS_i,\T_i)$
with maps $p_i:\T_i\rightarrow \T$
and $\bar\SS_i=\T_i\times_{\T} \bar\SS$,
such that for all $t\in \T$, and for all
$E_t\in \TS_{\SS_t}(A_t)$, there exists a $t'\in \T_i$ for some $i$
such that $(\bar\X_i)_{t'}\rightarrow (\bar\SS_i)_{t'}$ has generic fibre
$E_t$.

We will need five lemmas:

\proclaim Lemma 4.4. Let $i:\eta\rightarrow \SS$
be the inclusion of the generic
point of $\SS$ into $\SS$ in the situation
of Theorem 4.3, $g:\SS\rightarrow\T$ the restriction of $\bar g:
\bar \SS\rightarrow\T$. Then for any integer
$n>0$, ${}_nR^1g_*(i_*A)$, the $n$-torsion subsheaf of
$R^1g_*(i_*A)$, is a constructible sheaf on $\T$, where
$A$ is the
generic fibre of $f:\X\rightarrow\SS$.

Proof: Let $\Sigma$ be the discriminant locus of the elliptic
fibration $f:\X\rightarrow \SS$. Then
$$R^1g_*(i_*A)\subseteq R^1g_*'(i_*A)|_{\SS-\Sigma},$$
where $g'$ is the restriction of $g$ to $\SS-\Sigma$.
This can be verified on stalks since $\TS_S(A)\subseteq\TS_U(A)$ whenever
$U\subseteq S$ is an open subset of a variety $S$. Since a subsheaf of a
constructible
sheaf is constructible ([14, V 1.9]), it is enough to show the Lemma when
$\Sigma$ is empty, replacing $\SS$ by $\SS-\Sigma$.
Assume that $\Sigma$ is empty. By [3], 1.11 and 1.12,
we have
$$P_{\X/\SS}\cong i_*i^*P_{\X/\SS}$$
where $P_{X/S}=R^1f_*\Gm$ for a morphism $f:X\rightarrow S$,
and we also have the exact sequence
$$\exact{i_*A}{i_*i^*P_{\X/\SS}}{\boldz}.$$
([3], \S 1, (6))
Applying $g_*$, we see that
$$R^1g_*(i_*A)\cong R^1g_*(i_*i^*P_{\X/\SS})\cong R^1g_*P_{\X/\SS}.$$
Now by the Leray spectral sequence and [3], 1.4,
we have
$$R^1g_*P_{\X/\SS}\cong {R^2\pi_*\Gm\over R^2g_*\Gm},$$
where $\pi:\X\rightarrow\T$ is the restriction of
$\bar\pi:\bar\X\rightarrow\T$.
Taking $n$-torsion, we see that we obtain
an exact sequence
$$0\rightarrow {}_nR^2g_*\Gm\rightarrow {}_nR^2\pi_*\Gm
\rightarrow {}_nR^1g_*P_{\X/\SS}\rightarrow (R^2g_*\Gm)\otimes{\boldz/
n\boldz}\mapright{\phi} (R^2\pi_*\Gm)\otimes{\boldz/n\boldz},$$
and by the Kummer sequence we obtain
a diagram
$$\matrix{0&\mapright{}&(R^2g_*\Gm)\otimes{\boldz/n\boldz}&\mapright{}
&R^3g_*\mu_n\cr
&&\mapdown{\phi}&&\mapdown{f^*}\cr
0&\mapright{}&(R^2\pi_*\Gm)\otimes{\boldz/n\boldz}&\mapright{}
&R^3\pi_*\mu_n\cr}$$
Since $f^*$ is injective as $f$ has a section, $\phi$ is also
injective, and hence we obtain
$${}_nR^1g_*P_{\X/\SS}\cong {{}_nR^2\pi_*\Gm\over {}_nR^2g_*\Gm}.$$
Now we have by the Kummer sequence that
$R^2\pi_*\mu_n\rightarrow {}_nR^2\pi_*\G_m\rightarrow 0$,
and $R^2\pi_*\mu_n$ is constructible, so $(R^2\pi^*\Gm)_n$ is also.
Tracing back, we obtain the Lemma. $\bullet$

Let $f:X\rightarrow Y$ be a smooth, proper morphism, and let $D\subseteq X$
be an effective divisor, $D=D_1\cup\cdots\cup D_n$, with the $D_i$
irreducible components of $D$. We say $D$ is simple normal crossings
relative to $Y$ if whenever $x\in D$, and $I=\{i|x\in D_i\}$, then
$\bigcap_{i\in I} D_i\rightarrow Y$ is smooth of relative dimension
$\dim X-\dim Y-\#I$.

\proclaim Lemma 4.5. Let $f:X\rightarrow Y$ be a proper, smooth morphism
with $X$ and $Y$ of finite type over a field of characteristic zero. Let
$D\subseteq X$ be a divisor which has simple normal crossings over $Y$. Let
$\F$ be a locally constant constructible
torsion sheaf on $X$, and let $g:X-D\rightarrow Y$
be the restriction. Then $R^ig_*\F$ commutes with arbitrary base changes over
$Y$.

Proof. This is a standard application of purity (see e.g. [2],
Th. Finitude, Appendice).

First note that the lemma is true if $D$ is empty by the proper base change
theorem ([14], VI 2.3). Now let
$$X_r=X-\bigcup_{\# I=r}\left(\bigcap_{i\in I} D_i\right),$$
and $X_{n+1}=X$, where $n$ is the number of components of $D$. We then
have a diagram
$$\matrix{X_{r-1}&\hookrightarrow&X_r&\mapleft{i}&X_r-X_{r-1}=Z_r\cr
&\mapse{g}&\mapdown{f}&\mapsw{h}&\cr
&&S&&\cr}$$
Since $D$ is simple normal crossings over $Y$, $h$ is smooth.

By [14], VI 5.3, there is an exact sequence
$$\eqalign{\cdots\rightarrow R^{j-2c-1}h_*(i^*\F\otimes T_{Z_r/X_r})&
\rightarrow R^{j-1}f_*(\F|_{X_r})\rightarrow R^{j-1}g_*(\F|_{X_{r-1}})\cr
&\rightarrow R^{j-2c}h_*(i^*\F\otimes T_{Z_r/X_r})\rightarrow R^jf_*(\F|_{X_r})
\rightarrow\cdots\cr}$$
where $T_{Z_r/X_r}$ is a locally constant sheaf on $Z_r$. From this
we see that if $R^jg_*(\F|_{X_r})$
and $R^jh_*(i^*\F\otimes T_{Z_r/X_r})$ commute with base change for all $j$,
then so does $R^{j-1}g_*(\F|_{X_{r-1}})$. Thus by induction on $\dim X$ and
descending induction on $r$, we obtain the desired result. $\bullet$

\proclaim Lemma 4.6. Let $f:X\rightarrow Y$ be a compactifiable smooth
morphism, $X$ and $Y$ of finite type over an algebraically closed field
of characteristic zero. Then there exists a dense open set $U\subseteq Y$
such that for any locally constant constructible torsion sheaf $\F$ on
$X$, $(R^if_*\F)_{\bar y}=H^i(X_{\bar y},\F)$ for all $y\in U$, and
$(R^if_*\F)|_U$ is locally constant.

Proof. Let $X\subseteq \bar X\mapright{\bar f} Y$ be a compactification
of $f$. By resolution of singularities, we can assume that $\bar X-X$
is a simple normal crossings divisor. By generic smoothness, there exists
$U\subseteq Y$ such that $\bar X-X$ is simple normal crossings
over $U$. We then apply Lemma
4.5. The two claims then follow as in [14] VI 2.5 and 4.2. $\bullet$

Now consider the following situation. Let $f:X\rightarrow S$ be
an elliptic fibration with a section with generic fibre $A$, and let
$S'\subseteq S$ be an irreducible
locally closed subset such that the generic fibre $A'$
of $f':X\times_S S'\rightarrow S'$ is an elliptic curve.
Let $i$ be the inclusion
of the generic point of $S$ in $S$ and $i'$ the inclusion of the generic
point of $S'$ in $S'$. Let $g:S'\rightarrow S$ be the inclusion.
I claim there is a natural restriction map
$$\TS_S(A)\rightarrow\TS_{S'}(A')$$
defined as follows.
First, there is a map
$$i_*A\rightarrow g_*i'_*A'.$$
Indeed, a section of $i_*A$ corresponds to a rational section of
$f$ over an open set $U$. By restricting this rational section to
$U'=U\times_S S'$, we obtain a rational section of $f'$ over $U'$
provided that the original section was defined over an open dense set
of $S'$. But this is indeed the case, as a rational section is always
defined where $f$ is smooth by [20], Lemma 1.9.
Thus we obtain a composition $H^1(S,i_*A)\rightarrow H^1(S,g_*i_*'A')
\rightarrow H^1(S',i_*A')$, the latter map by the Leray spectral sequence.
This is the desired restriction map. It is then clear that
if $E\in\TS_S(A)$ and $g:Y\rightarrow S$ is an elliptic fibration with
generic fibre $E$, and if $g|_{S'}$ is an elliptic fibration  over $S'$,
then the latter
has generic fibre $E|_{S'}$, the image of $E$ in $\TS_{S'}(A')$.

As one particular case of this, there is, in our situation in Theorem 4.3,
a restriction map
$$(R^1g_*(i_*A))_{\bar t}\rightarrow \TS_{\SS_t}(A_t)$$
where $A_t$ is the generic fibre of $\X_t\rightarrow \SS_t$ whenever
$t\in\T$ is a closed point.
Indeed, $(R^1g_*(i_*A))_{\bar t}\cong H^1(\SS\times_{\T} \Spec\O_{\T,
\bar t}, i_*A)$ where $\O_{\T,\bar t}$ is the strict henselization of
the local ring of $\T$ at $t$, and we then use the above restriction map
with $S'=\SS_t$. We then have

\proclaim Lemma 4.7. In the situation of Theorem 4.3, there exists
a dense Zariski open subset $U\subseteq \T$
for which the following is true. For all closed points
$t\in U$, with $A_t$ the generic fibre of $\X_t\rightarrow \SS_t$,
the natural restriction map
$$({}_nR^1g_*(i_*A))_{\bar t}\rightarrow {}_n\TS_{\SS_t}(A_t)$$
is surjective, for any integer $n$.

Proof: First, let $\Sigma$ be the discriminant locus of $f$.
Let $g':\SS'=\SS-\Sigma\rightarrow\T$ be the restriction of $g$, $\X'=
\SS'\times_{\SS}\X$, and $f'$, $\pi'$ the restrictions of $f$ and
$\pi$ to $\X'$. As shown in the proof of Lemma 4.4, we have
$${}_nR^1g'_*(i_*A)|_{\SS'}\cong {{}_nR^2\pi_*'\Gm\over {}_nR^2g'_*\Gm}$$
and similarly
$${}_n\TS_{\SS'_t}(A_t)={{}_nH^2(\X'_t,\Gm)\over {}_nH^2(\SS'_t,\Gm)}.$$
To show the lemma for $\SS'$, it is then enough to show that there
is some open set $U\subseteq\T$ such that the natural
restriction map $${}_n(R^2\pi_*'\Gm)_{\bar t}\rightarrow
{}_nH^2(\X'_{\bar t},\Gm)$$ is surjective for $t\in U$ for any $n$.
But by the Kummer sequence, we have a diagram
$$\matrix{(R^2\pi_*'\mu_n)_{\bar t}&\mapright{}&{}_n(R^2\pi_*'\Gm)_{\bar t}&
\mapright{}&0\cr
\mapdown{}&&\mapdown{}&&\cr
H^2(\X'_{\bar t},\mu_n)&\mapright{}&{}_nH^2(\X'_{\bar t},\Gm)&
\mapright{}&0\cr}$$
By Lemma 4.5, the first vertical arrow is surjective for all $n$
for $t\in U$ for some
dense open set $U$ of $\T$, and hence so is the second. This proves the Lemma
for the family $(\X',\SS',\T)$.

Now note that
$$(R^1g'_*(i_*A)|_{\SS'})_{\bar t}=\TS_{\SS'(\bar t)}(A(\bar t))$$
where $\SS'(\bar t)=\SS'\times_{\T}\Spec\O_{\T,\bar t}$, and
$A(\bar t)=A\times_{\eta} \eta(\bar t)$, where $\eta(\bar t)$ is the generic
point of $\SS'(\bar t)$. We have a similar equality for $g$, and hence
a diagram
$$\matrix{0&\mapright{}&{}_n\TS_{\SS(\bar t)}(A(\bar t))&\mapright{}&
{}_n\TS_{\SS'(\bar t)}(A(\bar t))\cr
&&\mapdown{r_1}&&\mapdown{r_2}\cr
0&\mapright{}&{}_n\TS_{\SS_{\bar t}}(A_{\bar t})&\mapright{}&
{}_n\TS_{\SS'_{\bar t}}(A_{\bar t})\cr}$$
We have shown that
$r_2$ is surjective for $t\in U$. We wish to show that $r_1$ is.

Let $E\in {}_n\TS_{\SS_{\bar t}}(A_{\bar t})$,
$E'\in{}_n\TS_{\SS'(\bar t)}(A(\bar t))$
with $r_2(E')=E$. Let $Z\subseteq \SS(\bar t)$ be the closed subset of
$\SS(\bar t)$ defined by
$$Z=\{s\in\SS(\bar t)|loc_{\bar s}(E')\not=0\}.$$
Suppose $Z\not=\phi$. Let $s_1,\ldots,s_n$ be the generic points of the
irreducible components of $Z$. Since $\Sigma\rightarrow \T$ is smooth,
$\Sigma\times_{\T}\Spec \O_{\T,\bar t}$ in particular is non-singular,
and hence by [3], Cor. 3.2,  all the $s_i$ are generic points of
$\Sigma\times_{\T}\Spec\O_{\T,\bar t}$.
Thus $\overline{\{s_i\}}\cap\SS_{\bar t}\not=\phi$
(as each component of $\Sigma$ surjects onto $\T$),
and if
$s_i'$ is the generic point of this intersection, then $loc_{\bar s_i'}(E)
\not=0$,
contradicting $E\in\TS_{\SS_{\bar t}}(A_{\bar t})$. Thus $Z=\phi$
and $E'\in \TS_{\SS(\bar t)}(A(\bar t))$ and $r_1$ is surjective.
$\bullet$

\proclaim Lemma 4.8. With the hypothesis of Theorem 4.3, there exists
an $n\in\boldz$ and a dense open subset $U\subseteq\T$ such that
$\TS_{\SS_t}(A_t)$ is killed by $n$ for all closed points $t\in U$.

Proof: First we claim that there exists a dense open set $U\subseteq\T$
for which $$H^2(\X_{\bar t},\Gm)/H^2(\SS_{\bar t},\Gm)$$
is constant for all $t\in U$, assuming that it is always finite.
Recall from [8], II 3.2 that
if $X$ is a scheme smooth over a field $k$ of characteristic zero,
then $H^2(X,\Gm)=(\qz)^r\oplus G$
for some finite group $G$ and some $r$. We define $H^2(X,\Gm)^f:=G$. Thus since
we are assuming
$H^2(\X_{\bar t},\Gm)/H^2(\SS_{\bar t},\Gm)$ is finite,
this group is the same as
$H^2(\X_{\bar t},\Gm)^f/H^2(\SS_{\bar t},\Gm)^f$. Furthermore,
$H^2(X,\Gm)^f(l)=H^3(X,\boldz_l[1])_{tors}$ by [8], III, (8.9).
Thus by Lemma 4.5, there exists
an open subset $U\subseteq \T$ such that these groups are constant.

Now from the proof of [3], Theorem 2.24 a), there exists a locally constant
sheaf $\G$ on the discriminant locus $\Sigma$ of $f$ which is represented
by an \'etale scheme over $\Sigma$ such that there is an exact sequence
$$0\rightarrow H^2(\X_{\bar t},\Gm)/H^2(\SS_{\bar t},\Gm) \rightarrow
\TS_{\SS_{\bar t}}(A_t) \rightarrow H^1(\Sigma_{\bar t},\G|_{\SS_{\bar t}})$$
for all closed points $t\in \T$.
Again by Lemma 4.5, there is an open subset $U\subseteq \T$ on which the
latter group is constant (and finite). This proves the Lemma. $\bullet$

{\it Proof of Theorem 4.3:}
We prove the theorem by Noetherian induction, proving it for a dense
open subset $U\subseteq\T$. By [14, V, 1.8],
there exists an open subset $U\subseteq\T$ for which $\F={}_nR^1g_*(i_*A)|_U
$ is locally constant with finite stalks, (with $n$ given by Lemma 4.8),
as the latter sheaf is constructible by Lemma
4.4. Thus there exists an \'etale cover $U'\rightarrow U$ with
$\F|_{U'}$ constant, $\F|_{U'}\cong G_{U'}$, $G$ a finite group. Now
for all $x\in U',g\in G$, there exists an \'etale neighborhood $U_{x,g}
\rightarrow U'$
of $x$ such that the element of the stalk $g\in (\F|_{U'})_{\bar x}=G$ is
represented
by an element $E_{x,g}\in H^1(U_{x,g}\times_{\T} \SS,i_*A)$. By
quasi-compactness, we can find an open subcover of $\{U_{x,g}|x\in U'\}$
which covers $U'$,
for each $g$. In this manner, one obtains a finite collection of connected
\'etale  schemes $U_1,\ldots,U_n$ over $U'$ and elements
$E_j\in H^1(U_j\times_{\T}\SS,i_*A)$ such that for all closed points $t\in U$,
restricting each $E_j$ to $\SS_t$ gives all elements of $\TS_{\SS_t}(A_t)$,
by Lemmas 4.7 and 4.8.

We now only have to construct suitable models for the $E_i$. We will find a
decomposition of $U_i$ into a finite number of locally closed subsets $\T_{ij}$
along with families of elliptic fibrations
$(\bar\X_{ij},\bar\SS_{ij},\T_{ij})$,
$\bar\SS_{ij}\cong \T_{ij}\times_{\T}\bar\SS,$ such that the generic fibre of
$\bar\X_{ij}\rightarrow\bar\SS_{ij}$ is $E_i|_{\SS_{ij}}.$ This can be done
by constructing a family $(\X_{i1},\SS_{i1},\T_{i1})$ with
$\T_{i1}\subseteq U_i$ a dense open subset, and then applying Noetherian
induction, replacing $U_i$ with a finite number of non-singular locally
closed subsets of $U_i$ whose union is the complement of $\T_{i1}$ in
$U_i$, and then restricting $E_i$ to the inverse image of these locally
closed subsets on $\SS_i=\SS\times_{\T} U_i$.

One can construct a projective morphism
$f_i:\bar\X_i\rightarrow \bar\SS_i$, where $\bar\SS_i=\bar\SS\times_{\T} U_i$,
with generic fibre $E_i$, and by resolution of
singularities, $\bar\X_i$ can be taken to be non-singular. Now there is an
open set $\T_{i1}\subseteq U_i$ on which $\bar\X_i\rightarrow U_i$ is smooth,
by generic smoothness. This gives the desired family.
$\bullet$

{\it Proof of Theorem 0.1:} First apply Theorem 4.2 to a number of different
families, with $\L=\omega^{-1}_{\SS/\T}$:
\item{1)} $\T$ a point, $\SS$ a $\Ptwo$ or $F_e$, $0\le e\le 12$.
\item{2)} $\SS_i\rightarrow\T_i$ a finite set of families of surfaces which are
minimal resolutions of all possible rank one Gorenstein log Del Pezzo surfaces
with only $A_n$ singularities.
\item{3)} $\SS_i\rightarrow \T_i$ a finite set of families which include
all possible surfaces $S$ of the following type: $S$ is a minimal
resolution of a minimal ruled surface $S'\rightarrow \Pone$ with
$\le 4$ singular fibres with only $A_n$ singularities,
and there is a birational morphism
$S\rightarrow F_e$, $0\le e\le 2$.

It is not hard to construct the families in 2) and 3): in case 2), any
minimal resolution of a Gorenstein log Del Pezzo surface is the blowup
of $F_2$ in $\le 7$ points, by [17], Lemma 3. In this case, there is a
map $\SS_i\rightarrow\SS_i'$ over $\T_i$ induced by some sufficiently
high power of $\L$ which contracts only the $-2$ curves on the fibres
of $\SS_i\rightarrow\T_i$. Let $\E_i\subseteq\SS_i$ be the exceptional
locus of this contraction.

In case 3), any
such surface is a blowup of $F_e$, $0\le e\le 2$, in $\le 12$ points.
For example, to construct a family of resolutions which map to $F_e$
of minimal ruled surfaces
with two $A_1$ singularities and
one singular fibre, put $\T=F_e$, and obtain
$\SS$ by first blowing up $F_e\times\T$ along the diagonal $\Delta$ to
obtain $\alpha:\SS'\rightarrow F_e\times\T$ with exceptional locus
$\tilde\Delta$. Let $p:F_e\times\T\rightarrow\Pone\times\T$ be the
natural ruling. Next blow up $\SS'$ along $\tilde\Delta\cap F$, where $F$
is the proper transform of $p^{-1}(p(\Delta))$ via $\alpha$. This
new blow-up is $\SS$. Let $\E\subseteq\SS$ be the union of the
proper transforms of $\tilde\Delta$ and $F$ on $\SS$. This is the union of the
$-2$ curves on the fibres of $\SS\rightarrow\T$ which must be contracted
to obtain the singular, minimal ruled surfaces with one singular fibre
with two $A_1$ singularities. Likewise, one can construct other families
$\SS_i\rightarrow\T_i$ covering all such minimal ruled surfaces, and let
$\E_i$ be the corresponding union of $-2$ curves.

After applying Theorem 4.2, we obtain a new set of families $(\bar X_i,
\bar S_i,\T_i)$ in which in particular any jacobian Calabi-Yau elliptic
fibration appears. (We can throw out any family which contains threefolds
which aren't Calabi-Yau which may have popped up, as being Calabi-Yau
is deformation invariant.) One also has Weierstrass
models $w_i:\bar W_i\rightarrow \bar\SS_i$. Let $\Sigma_i\subseteq\bar\SS_i$
be the reduced discriminant locus of $\bar W_i$. By decomposing $\T_i$
into a finite number of locally closed subsets, we can assume that
if $\ZZ_i\subseteq\Sigma_i$ is the locus of points where $\Sigma_i\rightarrow
\T_i$ is not smooth, then each irreducible component of $\Sigma_i-\ZZ_i$
maps surjectively to $\T_i$. (For this to be the case, $\ZZ_i$ must
be of relative dimension zero.)
Furthermore, in cases 2) and 3) above, we have $\E_i\subseteq\SS_i$ which are
as above.
Now apply Theorem 4.3 with $\SS_i=\bar\SS_i-\ZZ_i\cup \E_i$.
By Theorem 3.5, the finiteness hypothesis of Theorem 4.3 is satisfied,
and Theorem 3.5 also tells us that all Calabi-Yau elliptic fibrations are
thus obtained in applying Theorem 4.3. $\bullet$

{\hd Bibliography}

\item{[1]} Deligne, P., ``Courbes Elliptiques: Formulaire d'apr\`es
J. Tate,'' in {\it Modular Functions in One Variable IV,} Lecture Notes
in Mathematics vol. 476, Springer-Verlag, 1975, pg. 53-74.
\item{[2]} Deligne, P., with Boutot, J.-F., Illusie, L., and Verdier,
J.-L., SGA $4{1\over 2}$, {\it Cohomologie \'Etale}, Lecture Notes
in Math., 569, Springer, Heidelberg 1977.
\item{[3]} Dolgachev, I., and Gross, M., ``Elliptic Three-folds I:
Ogg-Shafarevich Theory,'' preprint, 1992 (revised version, April 1993).
\item{[4]} Grassi, A., ``On Minimal Models of Elliptic Threefolds,''
{\it Math. Ann.} {\bf 290}, 287--301 (1991).
\item{[5]} Grassi, A., ``The Singularities of the Parameter Surface of
a Minimal Elliptic Threefold,'' To appear in {\it Int. J. of Math.,} {\bf 4},
???--??? (1993).
\item{[6]} Grassi, A., ``Log Contractions and Equidimensional Models of
Elliptic Threefolds,'' Preprint (1993).
\item{[7]} Gross, M., ``Elliptic Three-folds II: Multiple Fibres,'' MSRI
preprint \#018-93, 1992.
\item{[8]} Grothendieck, A., ``Le Groupe de Brauer, I, II, III,''
In {\it Dix Expos\'es
sur la Cohomologie des Sch\'emas,} North-Holland, Amsterdam, 1968, 46-188.
\item{[9]} Grothendieck, A., and Dieudonn\'e, J., {\it El\'ements de
G\'eom\'etrie Alg\'ebrique, I}, {\it Grundlehren der Math. Wissenschaften},
{\bf 166}, Springer-Verlag, (1971).
\item{[10]} Hunt, B., ``A Bound on the Euler Number for Certain
Calabi-Yau Threefolds,'' {\it J. Reine Angew. Math.,} {\bf 411} , 137--170
(1990).
\item{[11]} Kawamata, Y., ``Kodaira Dimensions of Certain Algebraic
Fiber Spaces,'' {\it J. Fac. Sci., Tokyo Univ. IA}, {\bf 30}, 1--24 (1983).
\item{[12]} Koll\'ar, J., and Mori, S., ``Classification of Three-dimensional
Flips,'' {\it J. of the AMS}, {\bf 5}, 533--703 (1992).
\item{[13]} Kollar, J., et al, {\it Flips and Abundance for
Algebraic Threefolds}, to appear.
\item{[14]} Milne, J. {\it \'Etale Cohomology}, Princeton Univ. Press, 1980.
\item{[15]} Miranda, R., ``Smooth Models for Elliptic Threefolds,'' in
{\it Birational Geometry of Degenerations,} Birkhauser, (1983) 85-133.
\item{[16]} Miranda, R., and Persson, U., ``On Extremal
Rational Elliptic Surfaces,'' {\it Math. Z.,} {\bf 193} 537--558 (1986).
\item{[17]} Miyanishi, M., and Zhang, D.Q., ``Gorenstein log Del Pezzo
Surfaces of Rank One,'' {\it Journal of Algebra,} {\bf 118}, 63--84,
(1988).
\item{[18]} Mumford, D., and Suominen, K., ``Introduction to the Theory of
Moduli,'' in {\it Algebraic Geometry, Oslo 1970,} Wolters-Noordhoff Press,
(1972) 171-222.
\item{[19]} Nakayama, N., ``On Weierstrass Models,'' {\it Algebraic
Geometry and Commutative Algebra in Honor of Masayoshi Nagata,} 405--431,
(1987).
\item{[20]} Nakayama, N., ``Local Structure of an Elliptic Fibration,''
Preprint, 1991, Univ. of Tokyo.
\item{[21]} Oguiso, K., ``On Algebraic Fiber Space Structures on a
Calabi-Yau Threefold,'' preprint, 1992.
\item{[22]} Sakai, F., ``The Structure of Normal Surfaces,'' {\it
Duke Math. J.,} {\bf 52}, 627--648 (1985).
\item{[23]} Schoen, C., ``On Fiber Products of Rational Elliptic
Surfaces with Section,'' {\it Math. Z.,} {\bf 197}, 177--199 (1988).
\end